\documentclass[11pt]{article}
\newcommand{\G}{\mathbb{G}}
\renewcommand{\baselinestretch}{1.2}
\newcommand{\dated}{\mbox{} \hfill {\small [{\tt \today}]}} \usepackage{amsmath,amssymb,amsfonts,diagrams}
\newenvironment{keywords}{\noindent\small {\it Keywords\/}:}{\vskip 4pt}
\newenvironment{classification}{\noindent\small 2000 {\it Mathematics Subject
Classification\/}:}{\vskip 12pt}

\newcommand{\comps}{{\mathbb C}}

\newcommand{\posints}{{\mathbb N}}

\newcommand{\tensor}{\otimes}

\newcommand{\cstar}{{C^\ast}}

\newcommand{\id}{{\mathrm{id}}}

\newcommand{\A}{{\mathfrak A}}

\newcommand{\Hilbert}{{\mathfrak H}}

\newcommand{\M}{{\mathfrak M}}

\newcommand{\VN}{\operatorname{VN}}

\newcommand{\varcl}[1]{\overline{#1}}

\usepackage{amsthm,enumerate}
\theoremstyle{plain}
\newtheorem{theorem}{Theorem}[section]
\newtheorem{lemma}[theorem]{Lemma}
\newtheorem{corollary}[theorem]{Corollary}
\newtheorem{proposition}[theorem]{Proposition}
\theoremstyle{definition}
\newtheorem{definition}[theorem]{Definition}
\theoremstyle{remark}

\newtheorem*{example}{Example}
\newtheorem*{rems}{Remarks}
\newtheorem*{exs}{Examples}

\newenvironment{examples}{\begin{exs}\begin{enumerate}}{\end{enumerate}\end{exs}}
\newenvironment{items}{\begin{enumerate}[\rm (i)]}{\end{enumerate}}
\newenvironment{alphitems}{\begin{enumerate}[\rm (a)]}{\end{enumerate}}

\title{Characterizations of compact and discrete quantum groups \\
through second duals}
\author{\textit{Volker Runde}\thanks{Research supported by NSERC under grant no.\ 227043-04.}}
\date{}
\begin{document}
\maketitle
\begin{abstract}
A locally compact group $G$ is compact if and only if $L^1(G)$ is an ideal in $L^1(G)^{\ast\ast}$, and the Fourier algebra $A(G)$ of
$G$ is an ideal in $A(G)^{\ast\ast}$ if and only if $G$ is discrete. On the other hand, $G$ is discrete if and only if ${\cal C}_0(G)$ is an ideal in ${\cal C}_0(G)^{\ast\ast}$. We show that these assertions are special cases of results on locally compact quantum groups in the sense of J.\ Kustermans and S.\ Vaes. In particular, a von Neumann algebraic quantum group $(\M,\Gamma)$ is compact if and only if $\M_\ast$ is an ideal in $\M^\ast$, and a (reduced) $\cstar$-algebraic quantum group $(\A,\Gamma)$ is discrete if and only if $\A$ is an ideal in $\A^{\ast\ast}$. 
\end{abstract}
\begin{keywords}
compact quantum group; discrete quantum group; locally compact quantum group; second dual;
weakly compact multiplication.
\end{keywords}
\begin{classification}
Primary 46L89; Secondary 22C05, 22D35, 43A99, 46H10, 46L51, 46L65, 47L50, 81R15, 81R50.
\end{classification}
\section*{Introduction}
Recently, J.\ Kustermans and S.\ Vaes introduced a surprisingly simple set of
axioms for what they call \emph{locally compact quantum groups}
(\cite{KV1} and \cite{KV2}). These axioms
cover both the Kac algebras (\cite{ES})---and thus all locally compact 
groups---as well as the compact quantum groups
in the sense of \cite{Wor} and allow for a notion of duality that extends the
Pontryagin duality for locally compact abelian groups.
\par
What makes the locally compact quantum groups interesting from the point of
view of abstract harmonic analysis is that many results in abstract harmonic
analysis can easily be
reformulated in the language of locally compact quantum groups. For instance,
Leptin's theorem (\cite{Lep}) asserts that a locally compact compact group
$G$ is amenable if and only if the Fourier algebra $A(G)$ (\cite{Eym})
has a bounded approximate identity. This statement can be rephrased as: $G$
is amenable if and only if its quantum group dual is co-amenable (as defined
in \cite[Definition 3.1]{BT}). Indeed, the ``if'' part is true for any
locally compact quantum group (\cite[Theorem 3.2]{BT}) whereas the converse
is known to be true only in the discrete case (\cite{Tom}; see also \cite{Rua} and \cite{BT}). 
\par
In \cite{Wat} (see also \cite{Gro} and \cite{Joh}), S.\ Watanabe proved that 
a locally compact group $G$ is compact if and only if its group algebra 
$L^1(G)$ is an ideal in its second dual. On the other hand, A.\ T.-M.\ Lau 
showed that $A(G)$ is an ideal in $A(G)^{\ast\ast}$ if and only if $G$ is 
discrete (\cite[Theorem 3.7]{Lau}). We shall see that both results are just 
two facets of one theorem on locally compact quantum groups. To this end, note that, according to
\cite{KV2}, one way of describing a locally compact quantum group is
as a Hopf--von Neumann algebra $(\M,\Gamma)$ with additional structure. Hence,
the predual $\M_\ast$ of $\M$ is a Banach algebra in a canonical way. 
Simultaneously extending the results by both Watanabe and Lau, we shall prove 
that the locally compact quantum group $(\M,\Gamma)$ is compact if and only if 
$\M_\ast$ is an ideal in its second dual.
\par 
If $G$ is a discrete group, then ${\cal C}_0(G) = c_0(G)$ is an ideal in ${\cal C}_0(G)^{\ast\ast} = \ell^\infty(G)$. On the other hand,
if ${\cal C}_0(G)$ is an ideal in ${\cal C}_0(G)^{\ast\ast}$, then multiplication in ${\cal C}_0(G)$ is weakly compact (\cite[Proposition 1.4.13]{Pal}) and thus compact (because ${\cal C}_0(G)$ has the Dunford--Pettis property), so that $G$ is discrete. We shall see that this result also extends to locally compact quantum groups: if $(\A,\Gamma)$ is a reduced $\cstar$-algebraic quantum group, then $(\A,\Gamma)$ is discrete if and only if $\A$ is an ideal in $\A^{\ast\ast}$.
\section{Compactness and amenability for Hopf--von Neumann algebras}
In this section, we formulate and briefly discuss the notions of compactness
and amenability in a general Hopf--von Neumann algebra context.
\par
We begin with recalling the definition of a Hopf--von Neumann algebra
($\bar{\tensor}$ denotes the $W^\ast$-tensor product):
\begin{definition}
A \emph{Hopf--von Neumann algebra} is a pair $(\M,\Gamma)$, where $\M$
is a von Neumann algebra and $\Gamma \!: \M \to \M \bar{\tensor} \M$ is
a \emph{co-multiplication}, i.e., a normal, unital $^\ast$-homomorphism
satisfying $(\id \tensor \Gamma) \circ \Gamma = (\Gamma \tensor \id) \circ
\Gamma$.
\end{definition}
\par
The main examples we have in mind are from abstract harmonic analysis:
\begin{examples}
\item For a locally compact group $G$, define $\Gamma_G \!: L^\infty(G) \to
L^\infty(G \times G)$ by letting
\[
  (\Gamma_G \phi)(x,y) := \phi(xy) \qquad (\phi \in L^\infty(G), \, x,y \in G).
\]
Then $(L^\infty(G),\Gamma_G)$ is a Hopf--von Neumann algebra.
\item For a locally compact group $G$, let $\lambda$ denote its left regular
representation, and let $\VN(G) := \lambda(G)''$ denote the
\emph{group von Neumann algebra} of $G$. Then 
\[
  \hat{\Gamma}_G \!: \VN(G) \to \VN(G) \bar{\tensor} \VN(G), \quad
  \lambda(x) \mapsto \lambda(x) \tensor \lambda(x)
\]
defines a co-multiplication, so that $(\VN(G), \hat{\Gamma}_G)$ is a 
Hopf--von Neumann algebra.
\end{examples}
\par
Given a Hopf--von Neumann algebra $(\M,\Gamma)$, one can define a product
$\ast$ on $\M_\ast$, the unique predual of $\M$, turning it into a Banach 
algebra:
\begin{equation} \label{prod}
  (f \ast g)(x) := (f \tensor g)(\Gamma x) \qquad (f,g \in \M_\ast, \,
  x \in \M).
\end{equation}
\begin{example}
Let $G$ be a locally compact group. Then applying (\ref{prod}) to
$(L^\infty(G), \Gamma_G)$ yields the usual convolution product on $L^1(G)$
whereas, for $(\VN(G), \hat{\Gamma}_G)$, we obtain pointwise multiplication 
on $A(G)$ (\cite{Eym}).
\end{example}
\par
For any Banach algebra $\A$, there are two canonical ways to extend the 
product to the second dual: the two \emph{Arens products} (see \cite{Dal}
or \cite{Pal}). These two products need not coincide; they are identical,
however, whenever one of the two factors involved is from $\A$. Given
a Hopf--von Neumann algebra $(\M,\Gamma)$, we will write $\ast$ for any of
the two Arens products on $\M_\ast^{\ast\ast} = \M^\ast$.
\par
We now define what it means for a Hopf--von Neumann algebra to be
(left) amenable and compact, respectively:
\begin{definition} \label{amdef}
Let $(\M,\Gamma)$ be a Hopf--von Neumann algebra. Then we call
$(\M,\Gamma)$ \emph{left amenable} if there is a \emph{left invariant}
state $M \in \M^\ast$, i.e., satisfying
\[
  f \ast M = f(1) M \qquad (f \in \M_\ast).
\]
If $M$ can be chosen to be in $\M_\ast$, we call $(\M,\Gamma)$
\emph{left compact}.
\end{definition}
\begin{example}
If $G$ is a locally compact group, then $(L^\infty(G),\Gamma_G)$ is
left amenable if and only if $G$ is amenable and left compact if and only if
$G$ is compact.
\end{example}
\par 
Of course one can equally well define 
right amenable (compact) Hopf--von Neumann 
algebras---through the existence of right invariant (normal) 
states---as well as amenable (compact) Hopf--von Neumann algebras by demanding
the existence of a (normal) state that is both left and right invariant.
\par
The following is well known, but for convenience, we include a proof:
\begin{proposition} \label{amprop}
Let $(\M,\Gamma)$ be a Hopf--von Neumann algebra. Then $(\M,\Gamma)$ is
amenable (compact) if and only if $(\M,\Gamma)$ is both left and right
amenable (compact).
\end{proposition}
\begin{proof}
If $M_l \in M^\ast$ is a left invariant state and $M_r \in \M^\ast$ is a
right invariant state, then $M_l \ast M_r$ is an invariant state (no matter
which Arens product is chosen on $\M^\ast$). Moreover, if $M_l$ and $M_r$ are
both normal, then so is $M_l \ast M_r$.
\end{proof}
\par
We conclude this section with recalling the notion of a co-involution.
\begin{definition}
Let $(\M,\Gamma)$ be a Hopf--von Neumann algebra. A \emph{co-involution}
for $(\M,\Gamma)$ is a $^\ast$-antihomomorphism $R \!: \M \to \M$ with
$R^2 = \id$ satisfying
\[
  (R \tensor R) \circ \Gamma = \sigma \circ \Gamma \circ R,
\]
where $\sigma$ is the flip map on $\M \bar{\tensor} \M$.
\end{definition}
\par
A co-involution $R$ for a Hopf--von Neumann algebra $(\M,\Gamma)$ is 
necessarily normal and can be used to define an involution $^\natural$
on the Banach algebra $\M_\ast$. For $f \in \M_\ast$, define $\bar{f}
\in \M_\ast$ by letting
\[
  \bar{f}(x) := \overline{f(x^\ast)} \qquad (x \in \M),
\] 
and set $f^\natural := \bar{f} \circ R$. It is easily seen that $^\natural$
is indeed an involution.
\par
In view of Proposition \ref{amprop}, we obtain:
\begin{corollary} \label{HvNcor}
The following are equivalent for a Hopf--von Neumann algebra $(\M,\Gamma)$
with co-involution:
\begin{items}
\item $(\M,\Gamma)$ is left amenable (compact);
\item $(\M,\Gamma)$ is right amenable (compact);
\item $(\M,\Gamma)$ is amenable (compact).
\end{items}
\end{corollary}
\section{Locally compact quantum groups}
Compact quantum groups  in the $\cstar$-algebraic setting were defined by
S.\ L.\ Woronowicz (see \cite{Wor} and \cite{MvD} for accounts). 
In \cite{KV1}, J.\ Kustermans and S.\ Vaes introduced a comparatively simple set of axioms
to define general locally compact quantum groups in a $\cstar$-algebraic
context. Alternatively, locally compact quantum groups can also be 
described as Hopf--von Neumann algebras with additional structure (see
\cite{KV2} and \cite{vDae}): both approaches are equivalent.
\par
In this section, we give an outline of the von Neumann algebraic approach to 
locally compact quantum groups. For details, see \cite{KV1}, \cite{KV2},
and also \cite{Kus}.
\par
We start with recalling some notions about weights on von Neumann algebras
(see \cite{Tak2}, for instance). 
\par
Let $\M$ be a von Neumann algebra, and let $\M^+$ denote its positive elements.
A \emph{weight} on $\M$ is an additive map $\phi \!: \M^+ \to [0,\infty]$ such
that $\phi(tx) = t \phi(x)$ for $t \in [0,\infty)$ and $x \in \M^+$. We let
\[
  {\cal M}_\phi^+ := \{ x \in \M^+ : \phi(x) < \infty \}, \qquad
  {\cal M}_\phi := \text{the linear span of ${\cal M}_\phi^+$},
\]
and
\[
  {\cal N}_\phi := \{ x \in \M : x^\ast x \in {\cal M}_\phi \}.
\]
Then $\phi$ extends to a linear map on ${\cal M}_\phi$, and ${\cal N}_\phi$
is a left ideal of $\M$. Using the GNS-construction (\cite[p.\ 42]{Tak2}),
we obtain a representation $\pi_\phi$ of $\M$ on some Hilbert space
$\Hilbert_\phi$; we denote the canonical map from ${\cal N}_\phi$ into
$\Hilbert_\phi$ by $\Lambda_\phi$. Moreover, we call $\phi$ 
\emph{finite} if ${\cal M}_\phi^+ = \M^+$, \emph{semi-finite}
if ${\cal M}_\phi$ is $w^\ast$-dense in $\M$, \emph{faithful} if $\phi(x) = 0$
for $x \in \M^+$ implies that $x = 0$, and \emph{normal} if $\sup_\alpha
\phi(x_\alpha) = \phi\left( \sup_\alpha x_\alpha \right)$ for each 
bounded, increasing net $( x_\alpha )_\alpha$ in $\M^+$. If $\phi$ is faithful
and normal, then the corresponding representation $\pi_\phi$ is faithful and
normal, too (\cite[Proposition VII.1.4]{Tak2}).
\begin{definition} \label{lcqg}
A \emph{locally compact quantum group} is a Hopf--von Neumann algebra
$(\M,\Gamma)$ such that:
\begin{alphitems}
\item there is a normal, semifinite, faithful 
weight $\phi$ on $\M$---a \emph{left Haar weight}---which is 
left invariant, i.e., satisfies
\[
  \phi((f \tensor \id)(\Gamma x)) = f(1) \phi(x) \qquad 
  (f \in \M_\ast, \, x \in {\cal M}_\phi);
\]
\item there is a normal, semifinite, faithful weight 
$\psi$ on $\M$---a \emph{right
Haar weight}---which is right invariant, i.e., satisfies
\[
  \phi((\id \tensor f)(\Gamma x)) = f(1) \psi(x) \qquad 
  (f \in \M_\ast, \, x \in {\cal M}_\psi).
\]
\end{alphitems}
\end{definition}
\begin{example}
Let $G$ be a locally compact group. Then the Hopf--von Neumann algebra
$(L^\infty(G),\Gamma_G)$ is a locally compact quantum group: $\phi$ and
$\psi$ can be chosen as left and right Haar measure, respectively.
\end{example}
\par
Even though only the existence of a left and a right Haar weight, respectively,
is presumed, both weights are actually unique up to a positive scalar multiple
(see \cite{KV1} and \cite{KV2}). 
\par
For each locally compact quantum group $(\M,\Gamma)$, there is a unique 
unitary---the \emph{multiplicative unitary}---$W \in 
{\cal B}(\Hilbert_\phi \tensor_2 \Hilbert_\phi)$, where $\tensor_2$
stands for the Hilbert space tensor product, such that
\[
  W^\ast(\Lambda_\phi(x) \tensor \Lambda_\phi(y))
  = (\Lambda_\phi \tensor \Lambda_\phi)((\Gamma y)(x \tensor 1))
  \qquad (x,y \in {\cal N}_\phi)
\]
(\cite[Theorem 1.2]{KV2}). The unitary $W$ lies in 
$\M \bar{\tensor} {\cal B}(\Hilbert_\phi)$
and implements the co-multiplication via 
\[
  \Gamma x = W^\ast (1 \tensor x) W \qquad (x \in \M)
\]
(see the discussion following \cite[Theorem 1.2]{KV2}). As 
discussed in \cite{KV1} and \cite{KV2}, locally compact quantum groups can 
equivalently be described in $\cstar$-algebraic terms. The
corresponding $\cstar$-algebra is obtained
as
\[
  \A := \varcl{\{ (\id \tensor \nu)(W) : 
                  \nu \in {\cal B}(\Hilbert_\phi)_\ast \}}^{\| \cdot \|}.
\]
\par
To emphasize the parallels between locally compact quantum groups and groups,
we shall use the following notation (which was suggested by Z.-J.\ Ruan): 
a locally compact quantum group $(\M,\Gamma)$ is denoted by the symbol
$\G$, and we write $L^\infty(\G)$ for $\M$, 
$L^1(\G)$ for $\M_\ast$, $L^2(\G)$ for
$\Hilbert_\phi$, and ${\cal C}_0(\G)$ for $\A$. If $L^\infty(\G) = L^\infty(G)$
for a locally compact group $G$ and $\Gamma = \Gamma_G$, we say that
$\G$ actually \emph{is} a locally compact group.
\par
The \emph{left regular} representation of a locally compact quantum group $\G$
is defined as
\[
  \lambda \!: L^1(\G) \to {\cal B}(L^2(\G)), 
  \quad f \mapsto (f \tensor \id)(W).
\]
(If $\G$ is a locally compact group $G$, this is just the usual left
regular representation of $L^1(G)$ on $L^2(G)$.)
It is a faithful, contractive representation of the Banach algebra $L^1(\G)$ on
$L^2(\G)$. 
\par
Locally compact quantum groups allow for the development of a duality
theory that extends Pontryagin duality for locally compact abelian groups.
\par
Set
\[
  L^\infty(\hat{\G}) : = 
  \varcl{\lambda(L^1(\G))}^{\,\text{$\sigma$-strongly$^\ast$}}.
\]
Then $L^\infty(\hat{\G})$ is a von Neumann algebra, and
\[
  \hat{\Gamma} \!:  L^\infty(\hat{\G}) \to 
   L^\infty(\hat{\G}) \bar{\tensor}  L^\infty(\hat{\G}), \quad
  x \mapsto \sigma W(x \tensor 1) W^\ast\sigma
\]
is a co-multiplication (here, $\sigma$ is the flip map on 
$L^2(\G) \tensor_2 L^2(\G)$).
One can also define a left Haar weight $\hat{\phi}$ and a right Haar weight 
$\hat{\psi}$ for $(L^\infty(\hat{\G}),\hat{\Gamma})$ turning it into a locally 
compact quantum group again, the \emph{dual quantum group} of $\G$, which
we denote by $\hat{\G}$. Generally, if
$X$ is an object associated with $\G$, we convene to denote the corresponding
object associated with $\hat{\G}$ by $\hat{X}$. Finally,
a Pontryagin duality theorem holds, i.e., $\Hat{\Hat{\G}} = \G$.
\begin{example}
If $\G$ is a locally compact group $G$, then $\hat{\G} =
(\VN(G), \hat{\Gamma}_G)$,
and $\hat{\phi} = \hat{\psi}$ is the Plancherel weight on $\VN(G)$
(\cite[Definition VII.3.2]{Tak2}).
\end{example}
\par
Each locally compact quantum group $\G$ is equipped with a canonical
co-involution---the \emph{unitary antipode}---, so that $L^1(\G)$ can be 
equipped with an involution $^\natural$. Unlike for groups, however,---and 
more generally for Kac algebras (see \cite{ES})---the left regular representation of 
$L^1(\G)$ need not be a $^\ast$-representation with respect to this involution.
\par
The \emph{antipode} of $\G$ is a $\sigma$-strongly$^\ast$-closed operator
$S$ on $L^\infty(\G)$ whose domain ${\cal D}(S)$ is 
$\sigma$-strongly$^\ast$-dense.
In general, $S$ is not totally defined. Letting
\[
  L^1_\ast(\G) := \{ f \in L^1(\G) : \text{there is $g \in L^1(\G)$
  such that $g(x) = \bar{f}(Sx)$ for $x \in {\cal D}(x)$} \},
\]
we obtain a dense subalgebra of $L^1(\G)$. On $L^1_\ast(\G)$, we can then
define an involution by letting
\[
  f^\ast(x) := \bar{f}(Sx) \qquad (x \in {\cal D}(S)).
\]
Defining a norm $||| \cdot |||$ on $L^1_\ast(\G)$ via
\[
  ||| f ||| := \max \{ \| f \|, \| f^\ast \| \} \qquad (f \in L^1_\ast(\G)),
\]
we turn $L^1_\ast(\G)$ into a Banach $^\ast$-algebra (in fact,
with isometric involution). The restriction
of $\lambda$ to $L^1_\ast(\G)$ is then a $^\ast$-representation of 
$L^1_\ast(\G)$. Furthermore, as a Banach $^\ast$-algebra, $L^1_\ast(\G)$
has an enveloping $\cstar$-algebra, which we denote by 
${\cal C}_0^u(\hat{\G})$.
\par
For more details, see \cite{Kus}
\section{Compact quantum groups}
We shall call a locally compact quantum group $\G$ \emph{compact} if
the Hopf--von Neumann algebra $(L^\infty(\G), \Gamma)$ is (left) compact in the
sense of Definition \ref{amdef}.
\par
Other characterizations of compactness for locally compact quantum groups
are collected below:
\begin{proposition} \label{cpchar}
The following are equivalent for a locally compact quantum group $\G$:
\begin{items}
\item $\G$ is compact;
\item the $\cstar$-algebra ${\cal C}_0(\G)$ is unital;
\item the left Haar weight of $\G$ is finite.
\end{items}
\end{proposition}
\begin{proof}
(i) $\Longleftrightarrow$ (ii) is \cite[Proposition 3.1]{BT} (compare
also \cite[Remark 3.8(5)]{Tom}).
\par
(ii) $\Longrightarrow$ (iii): Let $\phi$ denote the left Haar weight
of $\G$. Since ${\cal C}_0(\G) \cap {\cal N}_\phi$ is a norm dense left 
ideal of ${\cal C}_0(\G)$ is must contain the identity, so that 
$\phi$ is finite.
\par
(iii) $\Longrightarrow$ (i) is trivial.
\end{proof}
\par
In this section, we shall prove another,
less straightforward characterization of the compact quantum groups: a
locally compact quantum group $\G$ is compact if and only if $L^1(\G)$
is an ideal in $L^1(\G)^{\ast\ast}$.
\par
We first prove a lemma for general locally compact quantum groups:
\begin{lemma} \label{qlem}
Let $\G$ be a locally compact quantum group, and let $f \in L^1_\ast(\G)$.
Then we have
\[
  \lambda(f^\ast) \Lambda_\phi(x) = 
  \Lambda_\phi((\bar{f} \tensor \id)(\Gamma x)) \qquad (x \in {\cal N}_\phi).
\]
\end{lemma}
\begin{proof}
Fix $x \in {\cal N}_\phi$, and let $\nu \in {\cal B}(L^2(\G))_\ast$ be arbitrary.
Define $x\nu \in {\cal B}(L^2(\G))_\ast$ by letting
$(x\nu)(T) = \nu(Tx)$ for $T \in {\cal B}(L^2(\G))$. We obtain that
\[
  \begin{split}
  \nu(\lambda(f^\ast)x) & = (f^\ast \tensor \nu)(W(1 \tensor x)) \\
  & = (f^\ast \tensor x\nu)(W) \\
  & = \bar{f}( S(\id \tensor x\nu)(W)) ) \\
  & = \bar{f}( (\id \tensor x\nu)(W^\ast) ), \qquad
    \text{by \cite[Proposition 8.3]{KV1}}, \\
  & = (\bar{f} \tensor \nu)(W^\ast (1 \tensor x)) \\
  & = \nu((\bar{f} \tensor \id)(W^\ast)x).
  \end{split}
\]
Since $\nu \in {\cal B}(L^2(\G))_\ast$ is arbitrary, this means that
\[
  \lambda(f^\ast)x = (\bar{f} \tensor \id)(W^\ast)x.
\]
In view of the fact that $x \in {\cal N}_\phi$, it follows from
\cite[(8.3)]{KV1} that
\[
  \lambda(f^\ast)\Lambda_\phi(x) = (\bar{f} \tensor \id)(W^\ast)\Lambda_\phi(x)
  = \Lambda_\phi((\bar{f} \tensor \id) (\Gamma x)),
\]
as claimed.
\end{proof}
\par
Suppose now that $\G$ is a compact quantum group. Since the left Haar weight
$\phi$ is then finite by Proposition \ref{cpchar}---and can be chosen to 
be a state---, $\Lambda_\phi \!: {\cal N}_\phi \to L^2(\G)$ is a contraction
defined on all of $L^\infty(\G)$. Define $\iota \!:L^2(\G) \to L^1(\G)$ by
letting
\[
  \iota(\xi)(x) = \langle \xi, \Lambda_\phi(x^\ast) \rangle \qquad 
  (\xi \in L^2(\G), \, x \in L^\infty(\G)).
\]
Clearly, $\iota$ is a linear, injective contraction with dense range.
\par
The following compatibility result holds:
\begin{lemma} \label{cplem}
Let $\G$ be a compact quantum group. Then we have
\begin{equation} \label{compat}
  f \ast \iota(\xi) = \iota(\lambda(f)\xi) \qquad 
  (f \in L^1(\G), \, \xi \in L^2(\G)).
\end{equation}
\end{lemma}
\begin{proof}
By Lemma \ref{qlem} and the fact that $\lambda$ is a $^\ast$-representation
of $L^1_\ast(\G)$, we have
\begin{equation} \label{compeq}
  \lambda(f)^\ast \Lambda_\phi(x) = 
  \Lambda_\phi((\bar{f} \tensor \id)(\Gamma x)) \qquad 
  (f \in L^1_\ast(\G), \, x \in L^\infty(\G)).
\end{equation}
Hence, we thus obtain
\[
  \begin{split}
  \iota(\lambda(f) \xi)(x) & =  
  \langle \lambda(f)\xi, \Lambda_\phi(x^\ast) \rangle \\
  & = \langle \xi, \lambda(f)^\ast \Lambda_\phi(x^\ast) \rangle \\
  & = \langle \xi, \Lambda_\phi((\bar{f} \tensor \id)(\Gamma x^\ast)) \rangle,
      \qquad\text{by (\ref{compeq})}, \\
  & = \iota(\xi)(((\bar{f} \tensor \id)(\Gamma x^\ast))^\ast) \\
  & = \iota(\xi)((f  \tensor \id)(\Gamma x)) \\
  & = (f \tensor \iota(\xi))(\Gamma x) \\
  & = (f \ast \iota(\xi))(x) \qquad 
      (f \in L^1_\ast(\G), \, \xi \in L^2(\G), \, x \in L^\infty(\G)) 
  \end{split}
\]
and thus
\[
  \iota(\lambda(f) \xi) = f \ast \iota(\xi) \qquad 
  (f \in L^1_\ast(\G), \, \xi \in L^2(\G)).
\]
Since $L^1_\ast(\G)$ is dense in $L^1(\G)$, this yields (\ref{compat}).
\end{proof}
\par
Let $E$ and $F$ be Banach spaces, and let $T \!: E \to F$ be linear. Recall that $T$ is said to be \emph{weakly compact} if it maps the unital ball of $E$ onto a relatively weakly compact subset of $F$. 
\par 
As a consequence of Lemma \ref{cplem}, we obtain:
\begin{proposition} \label{wc}
Let $\G$ be a compact quantum group. Then multiplication in $L^1(\G)$
is weakly compact.
\end{proposition}
\begin{proof}
Let
\[
  I := \{ g \in L^1(\G) : \text{$L^1(\G) \ni f \mapsto f \ast g$ is weakly
  compact} \}. 
\]
We claim that $I = L^1(\G)$.
To see this, note first that---due to the reflexivity of $L^2(\G)$---the map
\[
  L^1(\G) \to L^2(\G), \quad f \mapsto \lambda(f) \xi
\]
is trivially weakly compact for each $\xi \in L^2(\G)$. Since 
$\iota \!: L^2(\G) \to L^1(\G)$ is continuous, Lemma 
\ref{cplem} yields that $\iota(L^2(\G)) \subset I$. Since
$I$ is closed and $\iota(L^2(\G))$ is dense in $L^1(\G)$, we conclude that
$I = L^1(\G)$.
\par
Since
\[
  f \ast g = (g^\natural \ast f^\natural)^\natural
  \qquad (f,g \in L^1(\G)),
\]
it follows that multiplication---both from the left and from the right---is
indeed weakly compact in $L^1(\G)$.
\end{proof}
\par 
It is a well known Banach algebraic fact that a Banach algebra is an ideal in its second dual if and only if multiplication in
$\A$ is weakly compact (\cite[Proposition 1.4.13]{Pal}). Hence, we obtain:
\begin{corollary} \label{compcor}
Let $\G$ be a compact quantum group. Then $L^1(\G)$ is an ideal in $L^1(\G)^{\ast\ast}$.
\end{corollary}
\par 
To prove the converse of Corollary \ref{compcor}, we first prove two lemmas, which are patterned on \cite[Lemmas 4.2 and 4.3]{MvD}. (We are grateful to S.\ Vaes for outlining this argument.)
\begin{lemma} \label{cplem1}
Let $\G$ be a locally compact quantum group such that $L^1(\G)$ is an ideal in $L^1(\G)^{\ast\ast}$. Then, for each state $f \in L^1(\G)$, there is a state $g \in L^1(\G)$ such that $f \ast g = g \ast f = g$.
\end{lemma}
\begin{proof}
Let $g \in L^1(\G)^{\ast\ast}$ be a $\sigma(L^1(\G)^{\ast\ast}, L^\infty(\G))$-accumulation point of the sequence $\left( \frac{1}{n} \sum_{k=1}^n f^{\ast k} \right)_{n=1}^\infty$, where $f^{\ast k}$ stands for the $k$-th power of $f$ with respect to the product $\ast$ in $L^1(\G)$. An inspection of the proof of \cite[Lemma 4.2]{MvD} yields that $g \in L^1(\G)^{\ast\ast}$ is a state satisfying $f \ast g = g \ast f = g$. Since $L^1(\G)$ is an ideal in $L^1(\G)^{\ast\ast}$, we see that $g \in L^1(\G)$.
\end{proof}
\begin{lemma} \label{cplem2}
Let $\G$ be a locally compact quantum group, and let $f,g \in L^1(G)$ be states such that $f \ast g = g \ast f = g$. Then, if $h \in L^1(\G)$ is such that $0 \leq h \leq f$, then $h \ast g = h(1) g$ holds.
\end{lemma}
\begin{proof}
The proof of \cite[Lemma 4.3]{MvD} carries over verbatim.
\end{proof}
\par
We can now state and prove the first main result of this paper:
\begin{theorem} \label{mainthm}
The following are equivalent for a locally compact quantum group $\G$:
\begin{items}
\item $\G$ is compact;
\item $L^1(\G)$ is an ideal in $L^1(\G)^{\ast\ast}$.
\end{items}
\end{theorem}
\begin{proof}
(i) $\Longrightarrow$ (ii) is Corollary \ref{compcor}.
\par
(ii) $\Longrightarrow$ (i): For each positive $f \in L^1(\G)$, set 
\[
  K_f := \{ g \in L^1(\G) : \text{$g$ is a state such that $f \ast g = f(1) g$} \}. 
\]
By Lemma \ref{cplem1}, $K_f$ is non-empty, and since multiplication with $f$ from the left is weakly compact, we see that $K_f$ is a weakly compact subset of $L^1(\G)$ whenever $f \neq 0$. Let $h,f \in L^1(\G)$ be states with $0 \leq h \leq f$. Then Lemma \ref{cplem2} yields that $K_f \subset K_h$. As a consequence, we have $K_{f_1+f_2} \subset K_{f_1} \cap K_{f_2}$ for any positive $f_1, f_2 \in L^1(\G)$. Fix a non-zero state $f_0 \in L^1(\G)$. Then $K_f \cap K_{f_0}$ is weakly compact for each positive $f \in L^1(\G)$ and non-empty (because it contains $K_{f+f_0}$). Consequently, $\bigcap \{ K_f : \text{$f \in L^1(\G)$ is positive} \}$ is non-empty. Any element of this intersection is a normal, left invariant state, so that $\G$ is compact.
\end{proof}
\par
Let $G$ be a locally compact group, and let $\G = (L^\infty(G),\Gamma_G)$.
Since $\G$ is compact if and only $G$ is compact, the main result of
\cite{Wat} is a particular case of Theorem \ref{mainthm}.
On the other hand, $\hat{\G} = (\VN(G),\hat{\Gamma}_G)$ is compact if and only 
if $G$ is discrete. From Theorem \ref{mainthm}, we thus 
recover \cite[Theorem 3.7]{Lau}:
\begin{corollary}
Let $G$ be a locally compact group. Then $A(G)$ is an ideal in 
$A(G)^{\ast\ast}$ if and only if $G$ is discrete.
\end{corollary}
\section{Discrete quantum groups}
There are various ways to define discrete quantum groups. Following
\cite{Tom}, we say that the locally compact quantum group $\G$ is
\emph{discrete} if $\hat{\G}$ is compact. As for compactness, we collect
a few equivalent characterizations:
\begin{proposition} \label{discprop}
The following are equivalent for a locally compact quantum group $\G$:
\begin{items}
\item $\G$ is discrete;
\item there is a central minimal projection in ${\cal C}_0(\G)$;
\item $L^1(\G)$ has an identity.
\end{items}
\end{proposition}
\begin{proof}
(i) $\Longrightarrow$ (ii) is noted in \cite[Remarks 3.8(2)]{Tom}.
\par
(ii) $\Longrightarrow$ (i): Let $p \in {\cal C}_0(\G)$ be a central 
minimal projection. Then $p$ is also central and minimal in $L^\infty(\G)$.
Let $\epsilon \in L^1(\G)$ denote the corresponding character. From the
definition of $\lambda$, it is clear that 
$\lambda(\epsilon) \in {\cal C}_0(\hat{\G})$ is a unitary operator on
$L^2(\hat{\G})$. Consequently, ${\cal C}_0(\hat{\G})$ has an identity, so that
$\hat{\G}$ is compact.
\par
(i) $\Longrightarrow$ (iii): In \cite[Remarks 3.8(2)]{Tom}, it is observed
that the character $\epsilon \in L^1(\G)$---as in (ii) $\Longrightarrow$ 
(i)---can also be chosen to satisfy
$(\epsilon \tensor \id) \circ \Gamma = (\id \tensor \epsilon) \circ \Gamma$
and thus is an identity for $L^1(\G)$.
\par
(iii) $\Longrightarrow$ (i): If $L^1(\G)$ has an identity, then so has
${\cal C}_0(\hat{\G}) = \overline{\lambda(L^1(\G))}$, so that $\hat{\G}$
is compact.
\end{proof}
\par
As a trivial consequence of Theorem \ref{mainthm}, a locally compact quantum group $\G$ is discrete if and only if $L^1(\hat{\G})$ is an ideal in $L^1(\hat{\G})^{\ast\ast}$. In this section, we give another characterization of discrete quantum groups in terms of second duals.
For the proof, we require a general fact on Banach $^\ast$-algebras, which may be of independent interest.
\par
We say that a $^\ast$-representation $\pi$ of a Banach $^\ast$-algebra $\A$ on a Hilbert space $\Hilbert$ is \emph{non-degenerate} if the closed linear span of $\{ \pi(a)\xi : a \in \A, \, \xi \in \Hilbert \}$ is all of $\Hilbert$.
\par
In our next lemma, let ${\cal K}(\Hilbert)$ denote the compact operators on a Hilbert space $\Hilbert$.
\begin{lemma} \label{Bstar}
Let $\A$ be a Banach $^\ast$-algebra, let $\Hilbert$ be a Hilbert space,
and let $\pi \!: \A \to {\cal B}(\Hilbert)$ be a faithful, non-degenerate
$^\ast$-representation of $\A$ such that $\pi(\A) \subset {\cal K}(\Hilbert)$.
Then, for each self-adjoint $a \in \A$, the spectra of $a$ in $\A$ and
of $\pi(a)$ in ${\cal B}(\Hilbert)$ coincide.
\end{lemma}
\begin{proof}
The claim is straightforward if $\dim \Hilbert < \infty$. We thus limit 
ourselves to the case where $\dim \Hilbert = \infty$.
\par
Let $a \in \A$ be self-adjoint, and let $\mathrm{Sp}(a)$ and 
$\mathrm{Sp}(\pi(a))$
denote the spectra of $a$ and $\pi(a)$ (in $\A$ and ${\cal B}(\Hilbert)$,
respectively). Since $\pi(a) \in {\cal K}(\Hilbert)$ and since $\dim 
\Hilbert = \infty$, we have $0 \in \mathrm{Sp}(\pi(a))$. Also, as $\A$ cannot
have an identity (again because $\pi(\A) \subset {\cal K}(\Hilbert)$),
we have $0 \in \mathrm{Sp}(a)$ as well. From \cite[Proposition 1.5.28]{Dal},
we conclude that $\mathrm{Sp}(\pi(a)) \subset \mathrm{Sp}(a)$.
\par
For the converse inclusion, let $\lambda \in \mathrm{Sp}(\pi(a))$, and suppose
without loss of generality that $\lambda \neq 0$. From the spectral theory
of compact operators on Hilbert space it is then well known that $\lambda$
is an eigenvalue of $\pi(a)$, i.e., there is $\xi \in \Hilbert \setminus
\{ 0 \}$ such that $(\lambda - \pi(a))\xi = 0$. Let $\A^\#$ denote
the \emph{unitization} of $\A$, i.e., the algebra $\A$ with an identity 
adjoined, and assume that $\lambda \notin \mathrm{Sp}(a)$. By the definition of
$\mathrm{Sp}(a)$ (see \cite[p.\ 78]{Dal}, for instance), there is $x \in \A^\#$
such that $x(\lambda -a) = 1$. The representation $\pi$ extends canonically
to a unital $^\ast$-representation $\pi^\# \!: \A^\# \to {\cal B}(\Hilbert)$,
so that we obtain 
\[
  \xi = \pi^\#(1)\xi = \pi^\#(x)(\lambda -\pi(a))\xi = 0,
\]
which is a contradiction.
\end{proof}
\par
By a $\cstar$-norm on a Banach $^\ast$-algebra $\A$ we mean a submultiplicative
norm $\gamma$ on $\A$ satisfying $\gamma(a^\ast a) = \gamma(a)^2$ for $a \in
\A$. We denote the completion of $\A$ with respect to $\gamma$ by 
$C^\ast_\gamma(\A)$, which is a $\cstar$-algebra. We say that $\A$ has
a \emph{unique $\cstar$-norm} if there is only one $\cstar$-norm on $\A$.
\par
We also require the notion of a regular Banach algebra (\cite[Definition 4.1.16]{Dal}): a commutative Banach algebra $\A$ with charater space $\Phi_\A$ is called \emph{regular} if, for each closed $F \subset \Phi_\A$ and $\phi \in \Phi_\A \setminus F$, there is $a \in \A$
with $\hat{a} |_F \equiv 0$ and $\hat{a} (\phi) = 1$, where $\hat{a}$ stands for the Gelfand transform of $a$.
\begin{proposition} \label{Bprop}
Let $\A$ be a Banach $^\ast$-algebra, and suppose that there is a $\cstar$-norm
$\gamma$ on $\A$ such that $C^\ast_\gamma(\A)$ is an ideal in 
$C^\ast_\gamma(\A)^{\ast\ast}$. Then $\A$ has a unique $\cstar$-norm.
\end{proposition}
\begin{proof}
Suppose without loss of generality that $\A$ is infinite-dimensional.
\par
Let $a \in \A$ be self-adjoint, and let $\langle a \rangle$ denote the
closed subalgebra of $\A$ generated by $a$. We claim that the Banach algebra $\langle a \rangle$
is regular.
\par
Since $C^\ast_\gamma(\A)$ is an ideal in $C^\ast_\gamma(\A)^{\ast\ast}$, 
it follows
from \cite[Exercises III.5.3 and III.5.4]{Tak}, that there is a family
$( \Hilbert_\alpha )_\alpha$ of Hilbert spaces such that
\[
  C^\ast_\gamma(\A) 
  \cong \text{$c_0$-}\bigoplus_\alpha {\cal K}(\Hilbert_\alpha)
  \subset {\cal K}(\Hilbert),
\]
where $\Hilbert := \text{$\ell^2$-}\bigoplus_\alpha \Hilbert_\alpha$. This
yields a faithful, non-degenerate $^\ast$-representation $\pi \!: \A \to
{\cal B}(\Hilbert)$ with $\pi(\A) \subset {\cal K}(\Hilbert)$. From Lemma
\ref{Bstar} and the spectral theory of compact operators on Hilbert space
we conclude that the spectrum of $a$ in $\A$ is of the form $\{ \lambda_n :
n \in \posints \} \cup \{ 0 \}$, where $\lambda_n \neq 0$ for $n \in \posints$
and $\lim_{n \to \infty} \lambda_n = 0$. Consequently, the character space
of $\langle a \rangle$ is canonically homeomorphic to the discrete subset
$\{ \lambda_n : n \in \posints \}$ of $\comps$. With the help of the 
holomorphic functional calculus, it is then straightforward to see that 
$\langle a \rangle$ is regular, as claimed.
\par
Since $a \in \A$ was an arbitrary self-adjoint element, $\A$ is 
therefore \emph{locally regular} in the sense of \cite[Definition 4.1]{Bar}
and thus has a unique $\cstar$-norm (\cite[Lemma 4.2]{Bar}).
\end{proof}

\begin{theorem} \label{mainthm2}
The following are equivalent for a locally compact quantum group $\G$:
\begin{items}
\item $\G$ is discrete;
\item ${\cal C}_0(\G)$ is an ideal in ${\cal C}_0(\G)^{\ast\ast}$.
\end{items}
\end{theorem}
\begin{proof}
(i) $\Longrightarrow$ (ii): Let 
\[
  I := \{ a \in {\cal C}_0(\G): \text{${\cal C}_0(\G) \ni b \mapsto ba$ is weakly
  compact} \},
\] 
and note that $I$ is closed. Since $\G$ is discrete, $\hat{\G}$ is compact, so that we have the map $\iota \!: L^2(\hat{\G}) \to L^1(\hat{\G})$ as defined immediately prior to Lemma \ref{cplem}; moreover, we have
\[
  \hat{\lambda}(f) \hat{\lambda}(\iota(\xi)) = \hat{\lambda}(f \ast \iota(\xi)) = \hat{\lambda}(\iota(\hat{\lambda}(f)\xi))
  \qquad (f \in L^1(\hat{\G}), \, \xi \in L^2(\hat{\G})).
\]
Since $L^2(\hat{\G})$ is reflexive, it follows---as in the proof of Proposition \ref{wc}---that $\hat{\lambda}(\iota(L^2(\hat{\G})))$ is contained in $I$. Since $\hat{\lambda}(\iota(L^2(\hat{\G})))$ is dense in ${\cal C}_0(\G)$, $I$ must equal all of ${\cal C}_0(\G)$, i.e., right multiplication in ${\cal C}_0(\G)$ is weakly compact. Using the involution on ${\cal C}_0(\G)$, we see that left multiplication in ${\cal C}_0(\G)$ is also weakly compact, so that ${\cal C}_0(\G)$ is an ideal in ${\cal C}_0(\G)^{\ast\ast}$.
\par
(ii) $\Longrightarrow$ (i): Let $\gamma$ denote the $\cstar$-norm on $L^1_\ast(\hat{\G})$ induced by $\hat{\lambda}$, so that
$C^\ast_\gamma(L^1_\ast(\hat{\G})) = {\cal C}_0(\G)$. From Proposition \ref{Bprop}, we conclude that $L^1_\ast(\hat{\G})$ has a unique $\cstar$-norm; in particular, ${\cal C}_0^u(\G) \cong {\cal C}_0(\G)$ holds. Consequently, ${\cal C}_0(\G)$ has a character, say $\epsilon$, by \cite[Theorem 3.1]{BT}. Let $p \in {\cal C}_0(\G)^{\ast\ast}$ be the corresponding central minimal projection. Choose $a \in {\cal C}_0(\G)$ with $\epsilon(a) = 1$. It follows that $p = \epsilon(a) p = ap \in {\cal C}_0(\G)$. By Proposition \ref{discprop}, this means that $\G$ is discrete.
\end{proof}
\par 
Since compactness and discreteness are dual to each other, we finally obtain as a consequence of both Theorems \ref{mainthm} and \ref{mainthm2}:
\begin{corollary} \label{lastcor}
The following are equivalent for a locally compact quantum group $\G$:
\begin{items}
\item $\G$ is compact;
\item $L^1(\G)$ is an ideal in $L^1(\G)^{\ast\ast}$;
\item ${\cal C}_0(\hat{\G})$ is an ideal in ${\cal C}_0(\hat{\G})^{\ast\ast}$;
\item $\hat{\G}$ is discrete.
\end{items}
\end{corollary}
\renewcommand{\baselinestretch}{1.0}
\renewcommand{\baselinestretch}{1.2}
\dated
\vfill
\begin{tabbing}
{\it Author's address\/}: \= Department of Mathematical and Statistical Sciences \\
\> University of Alberta \\
\> Edmonton, Alberta \\
\> Canada T6G 2G1 \\[\medskipamount]
{\it E-mail\/}: \> {\tt vrunde@ualberta.ca} \\[\medskipamount]
{\it URL\/}: \> {\tt http://www.math.ualberta.ca/$^\sim$runde/} 
\end{tabbing} 
\end{document}